\documentclass[12pt]{article}
\usepackage[latin1]{inputenc}
\usepackage{amssymb,amsmath,amsfonts}
\newtheorem{theorem}{Theorem}
\newtheorem{proposition}[theorem]{Proposition}
\newtheorem{lemma}[theorem]{Lemma}

\newtheorem{corollary}[theorem]{Corollary}
\newtheorem{fact}[theorem]{Fact}

\newtheorem{remarks}[theorem]{Remarks}

\newcommand{\R}{\mathbb{R}}
\newcommand{\U}{\mathcal{U}}

\newcommand{\C}{\mathbb{C}}

\newcommand{\spa}{\mbox{span\,}}

\newcommand{\rank}{\mbox{rank }}
\newcommand{\kerl}{\mbox{ker}}

\newcommand{\po}{{\hspace*{-1ex}}{\bf .  }}

\def\lp{{\langle\!\langle}}
\def\rp{{\rangle\!\rangle}}
\def\<{{\langle}}
\def\>{{\rangle}}

\def\Tal{{\cal T}}

\def\a{\alpha}
\def\be{\begin{equation} }
\def\ee{\end{equation} }

\def\proof{\noindent\emph{Proof: }}
\def\qed{\ifhmode\unskip\nobreak\fi\ifmmode\ifinner
\else\hskip5 pt \fi\fi\hbox{\hskip5 pt \vrule width4 pt
height6 pt  depth1.5 pt \hskip 1pt }}

\begin{document}

\title{Holomorphicity of real Kaehler submanifolds}
\author{A. de Carvalho, S. Chion and M. Dajczer}
\date{}
\maketitle

\begin{abstract} 
Let $f\colon M^{2n}\to\R^{2n+p}$ denote an isometric immersion
of a Kaehler manifold of complex dimension $n\geq 2$ into 
Euclidean space with codimension $p$. 
If \mbox{$2p\leq 2n-1$}, we show that generic rank conditions on the 
second fundamental form of the submanifold imply that $f$ has to 
be a minimal submanifold. In fact, for codimension $p\leq 11$ 
we prove that $f$ must be holomorphic with respect to some 
complex structure in the ambient space.  
\end{abstract}

\renewcommand{\thefootnote}{\fnsymbol{footnote}} 
\renewcommand{\thefootnote}{\arabic{footnote}} 

\renewcommand{\thefootnote}{\fnsymbol{footnote}} 
\renewcommand{\thefootnote}{\arabic{footnote}}  

Throughout the paper $f\colon M^{2n}\to\R^{2n+p}$ denotes a \emph{real 
Kaehler submanifold}. This means that $(M^{2n},J)$ is a connected Kaehler 
manifold of complex dimension $n\geq 2$  isometrically immersed into 
Euclidean space with substantial codimension $p$. The latter condition says 
that the codimension cannot be reduced, even locally.  It is well-known 
that real Kaehler submanifolds in low codimension are generically holomorphic. 
This means that $p$ is even and $f$ is holomorphic with respect to some 
complex structure in $\R^{2n+p}$. In fact, it was shown by Dajczer and 
Rodr\'{i}guez \cite{DR1} that this is the case if the type number of the 
immersion is $\tau(x)\geq 3$ at any point $x\in M^{2n}$.  Notice that 
this strong assumption forces the codimension of the immersion to satisfy 
$3p\leq 2n$.

In view of the above, it is a natural task to understand what can (locally) 
happen in low codimension under weaker, and algebraically simpler, assumptions 
on the second fundamental form than the one on the type number. In view 
of the rigidity results given in \cite{CD}, it seems natural to hope that 
generic conditions imply holomorphicity just for codimension 
$2p\leq 2n-1$. 
Moreover, it is desirable to replace the type number assumption for something 
simpler and more meaningful. By the latter, we mean providing a workable
condition as a starting point in order to study the real Kaehler submanifold
that are not holomorphic.

In fact, for very low codimensions there is already some relevant 
knowledge in the direction pointed above.  For instance, Dajczer 
\cite{Da} showed that if the index of relative nullity of 
$f\colon M^{2n}\to\R^{2n+2}$ satisfies \mbox{$\nu(x)<2(n-2)$} at any 
point $x\in M^{2n}$, then the submanifold is holomorphic along each 
connected component of an open dense subset of $M^{2n}$.
Recall that the \emph{index of relative nullity} $\nu(x)$ of $f$ at 
$x\in M^{2n}$ is the dimension of the relative nullity tangent subspace
defined by
$$
\Delta(x)=\mathcal{N}(\a)(x)=\{X\in T_xM:\a(X,Y)
=0\;\;\mbox{for all}\;\;Y\in T_xM\},
$$
where $\a\colon TM\times TM\to N_fM$ stands for the second fundamental 
form of $f$.  The study of non-holomorphic real Kaehler submanifolds 
$f\colon M^{2n}\to\R^{2n+2}$ with index of relative nullity 
$\nu(x)=2(n-2)$ was done in \cite{DG2}, \cite{DG3} and \cite{DR2}.

If $f\colon M^{2n}\to\R^{2n+3}$ satisfies $\nu(x)<2(n-3)$ at any  
$x\in M^{2n}$,  Dajczer and Gromoll \cite{DG4}  proved that  there 
exists an open dense subset $U$ of $M^{2n}$ 
such that, along each connected component  $U'$ of $U$, the submanifold 
$f|_{U'}$ has a Kaehler extension, namely, there exist a real Kaehler 
hypersurface $j\colon N^{2n+2}\to\R^{2n+3}$ and a holomorphic isometric 
immersion $h\colon U'\to N^{2n+2}$ such that $f|_{U'}=j\circ h$.
We point out that real Kaehler hypersurfaces have been classified 
by Dajczer and Gromoll \cite{DG1} by means of the Gauss parametrization
in terms of pseudoholomorphic surfaces in spheres.

Yan and Zheng \cite{YZ} observed that both results discussed
above still hold under the slightly weaker assumption that the 
\emph{complex index of relative nullity}, defined by
$$
\nu^c(x)=\dim\Delta(x)\cap J\Delta(x),
$$ 
satisfies the same pointwise inequalities required for $\nu(x)$. The 
main result in  \cite{YZ} is that if $f\colon M^{2n}\to\R^{2n+4}$ satisfies 
$\nu^c(x)<2(n-4)$ at any $x\in M^{2n}$, then there exists an open dense 
subset $U$ of $M^{2n}$ such that along each connected component $U'$ of $U$, 
the submanifold $f|_{U'}$ has a Kaehler extension, namely, there exist a 
real Kaehler submanifold $j\colon N^{2n+2}\to\R^{2n+4}$ 
and a holomorphic isometric immersion $h\colon U'\to N^{2n+2}$ such that 
$f|_{U'}=j\circ h$. Moreover, although the extension $j$ may not be unique, 
it can be chosen to be minimal if $f$ is minimal. Of course, we may be in
the situation where $f$ itself is holomorphic. Finally, the case of 
codimension $p=6$ was considered by Carvalho and Guimar\~aes \cite{CG}.

Since its introduction by do Carmo and Dajczer \cite{CD}, the notion 
of $s$-nullity of an isometric immersion has played a leading role in 
the understanding of rigidity questions of submanifolds in low 
codimension.  A complex version for real Kaehler submanifolds of this 
notion goes as follows.
\vspace{1ex}

The \emph{complex $s$-nullity} $\nu^c_s(x)$ of $f\colon M^{2n}\to\R^{2n+p}$
at $x\in M^{2n}$, $1\leq s\leq p$, is 
$$
\nu^c_s(x)=\max_{U^s\subset N_fM(x)}\{\dim(\mathcal{N}
(\a_{U^s})\cap J\mathcal{N}(\a_{U^s}))\},
$$
where $U^s$ denotes an $s$-dimensional subspace of $N_fM(x)$
and $\a_{U^s}=\pi_{U^s}\circ\a$ being $\pi_{U^s}\colon N_fM\to U^s$
the projection. Notice that $\nu^c_p(x)=\nu^c(x)$.

\begin{theorem}\po\label{minimal} Let $f\colon M^{2n}\to\R^{2n+p}$ 
be a real Kaehler submanifold.  Assume that the complex $s$-nullities  
satisfy $\nu^c_s(x)<2(n-s)$ at any $x\in M^{2n}$ for all $1\leq s\leq p$. 
Then $f$ is a minimal immersion.  
\end{theorem}

Minimal real Kaehler submanifolds have been intensively studied since
it was shown in \cite{DG4} that they posses several of the basic 
properties of minimal surfaces. For instance, any simply connected 
minimal real Kaehler submanifold $f\colon M^{2n}\to\R^{2n+p}$ is either 
holomorphic or has a non trivial one-parameter associated family of 
minimal isometric immersions with the same Gauss map. Moreover, $f$ 
can be realized as the ``real part'' of its holomorphic representative 
$F\colon M^{2n}\to\C^{2n+p}$, namely, $\sqrt{2}F=(f,\bar{f})$ where 
$\bar{f}\colon M^{2n}\to\R^{2n+p}$ is the conjugate immersion to $f$ 
in the associated family.
\vspace{1ex}

Notice that already for codimension $p=4$ the composition of isometric 
immersions obtained in the aforementioned result by Yan and Zheng shows 
that Theorem \ref{minimal} does not hold if we drop just one of the 
assumptions on the complex $s$-nullities.
\vspace{1ex}

Once one obtains minimality from Theorem \ref{minimal}, one can 
make use of the rigidity theorem for isometric immersions of 
codimension $p\leq 5$ given in \cite{CD} to conclude that the submanifold 
must, in fact, be holomorphic; see Remarks \ref{remark} for details. 

\begin{theorem}\label{main}\po
Let $f\colon M^{2n}\to\R^{2n+p}$ be a real Kaehler submanifold. Assume 
that the codimension is $p\leq 11$ and that the complex $s$-nullities  
satisfy $\nu^c_s(x)<2(n-s)$ at any $x\in M^{2n}$ for all $1\leq s\leq p$. 
Then $p$ is even and $f$ is holomorphic.  	
\end{theorem}

The examples of minimal but not holomorphic real Kaehler submanifold with
codimension $p=2$ constructed in \cite{DG2} satisfy one of the assumptions 
on the complex $s$-nullities in the above result since $\nu^c_1(x)=2n-4$
but not the other one since $\nu^c_2(x)=2n-4$.
\vspace{1ex}

It remains an open problem if Theorem \ref{main} still holds if 
the assumption $p\leq 11$ is dropped. In that respect, we observe 
that the core of the proof of Theorem \ref{main} is a result in the 
theory of flat bilinear forms that only holds until dimension $11$.
In fact, we conclude the paper constructing examples showing
that this result is false if the dimension is $12$.
On one hand, if  Theorem \ref{main} holds for any codimension,
then, of course,  Theorem \ref{minimal} follows from 
Theorem \ref{main}. On the other hand, Theorem \ref{minimal}
is a fundamental ingredient in our proof of Theorem \ref{main}.

\section{The proofs}

The proofs of our theorems rely heavily on results in the 
realm of the theory of flat bilinear forms tailored for that 
purpose. Hence, we first recall from that theory as lemmas 
several facts already in the literature.
\vspace{1ex}

Let $V^n$ and $W^m$ denote (real) vector spaces of finite 
dimension $n$ and $m$, respectively. Given a (maybe not symmetric) 
bilinear form $\beta\colon V^n\times V^n\to W^m$, we denote by
$$
\mathcal{S}(\beta)=\spa\{\beta(X,Y)\colon X,Y\in V^n\}
$$
the subspace generated by $\beta$ and its (right) kernel by
$$
\mathcal{N}(\beta)=\{Y\in V^n\colon\beta(X,Y)=0
\;\mbox{for all}\;X\in V^n\}.
$$ 
We also denote
\be\label{r}
r=\max\{\dim B_X(V)\colon X\in V^n\},
\ee
where $B_X\colon V^n\to W^m$ is the linear transformation 
defined by $B_XY=\beta(X,Y)$.  
A vector $X\in V^n$ is called a (left) \emph{regular element} 
of $\beta$ if $\dim B_X(V)=r$. The set $RE(\beta)$ of regular 
elements of $\beta$ is easily seen to be an open dense subset 
of $V^n$; for instance see  Proposition $4.4$ in \cite{DT}. 

\begin{lemma}\po\label{l:symmetricspan}  If 
$\beta\colon V^n\times V^n\to W^m$ is symmetric and 
$\mathcal{S}(\beta)=W^m$ then  $r(r+1)\geq 2m$.  
\end{lemma}

\proof On one hand, the set of non-asymptotic regular 
elements
$$
RE^*(\beta)=\{X\in RE(\beta):\beta(X,X)\neq 0\}
$$
is open and dense in $V^n$. On the other hand, from Eq.\ (8)
in \cite{Mo} or Proposition $4.6$ in \cite{DT} we have
$\beta(V^n,Z)\subset B_X(V^n)$ for any $Z\in\kerl B_X$.
Hence, given $Y_1\in RE^*(\beta)$ consider $Y_2,\ldots,Y_r$ 
such that
$$
B_{Y_1}(V^n)=\spa\{\beta(Y_1,Y_j),\;1\leq j\leq r\}.
$$
Then 
$$
W^m=\spa\{\beta(Y_i,Y_j),\;1\leq i,j\leq r\},
$$
and the result follows.\vspace{2ex}\qed

Let $W^{p,p}$ denote a $2p$-dimensional real vector space endowed 
with an inner product of signature $(p,p)$. Hence $p$ is the dimension 
of the vector subspaces of $W^{p,p}$ of maximal dimension where 
the induced inner product is either positive or negative definite. 
A vector subspace $L\subset W^{p,p}$ is called \emph{degenerate} if 
$\U=L\cap L^\perp\neq\{0\}$ and \emph{nondegenerate} otherwise. 
We call the \emph{rank} of $L$ the rank of the inner 
product induced on $L$, that is, $\rank L=\dim L-\dim\U$.
If $L=\U$, then $L$ is called an \emph{isotropic} subspace. 

\begin{lemma}\po
Given a subspace $L\subset W^{p,p}$ there is a direct sum decomposition 
\be\label{e:decomposition}
W^{p,p}=\U\oplus\hat{\U}\oplus\mathcal{V}
\ee
such that $\U=L\cap L^\perp$, the subspace $\hat{\U}$ is isotropic 
with $L\subset\U\oplus\mathcal{V}$ and
$\mathcal{V}=(\U\oplus\hat{\U})^\perp$ is a nondegenerate subspace.
\end{lemma}

\proof See Sublemma $2.3$ in \cite{CD} or Corollary 4.3 in \cite{DT}.
\qed

\begin{lemma}\po\label{p:densitydeg} 
Given $\beta\colon V^n\times V^n\to W^{p,p}$ set  
$$
RE^o(\beta)=\{X\in RE(\beta)\colon\dim\U(X)=\tau\},
$$
where $\U(X)=B_X(V)\cap B_X(V)^\perp$ and 
$$
\tau=\min\{\dim\U(X)\colon X\in RE(\beta)\}.
$$ 
Then $RE^o(\beta)$ is an open dense subset of $V^n$.
\end{lemma}
 
\proof See the proof of Lemma 2.1 in \cite{DR1}.\vspace{2ex}\qed
\newpage

A bilinear form $\beta\colon V^n\times V^n\to W^{p,p}$ is called 
\emph{flat} if
$$
\<\beta(X,Y),\beta(Z,T)\>-\<\beta(X,T),\beta(Z,Y)\>=0
$$
for all $X,Y,Z,T\in V^n$. It is said that $\beta$ is \emph{null} when
$$
\<\beta(X,Y),\beta(Z,T)\>=0
$$
for all $X,Y,Z,T\in V^n$. Thus null bilinear forms are flat.
\vspace{1ex}

\begin{lemma}\po\label{bilfla}
Let $\beta\colon V^n\times V^n\to W^{p,p}$ be a flat bilinear form. 
If $X\in RE(\beta)$ then 
\be\label{bili}
\spa\{\beta(Y,Z): Y\in V^n\;\;\text{and}\;\;Z\in\ker B_X\}\subset\U(X).
\ee 
\end{lemma}

\proof See Sublemma $2.4$ in \cite{CD} or Proposition $4.6$ in 
\cite{DT}.\qed

\begin{proposition}\po\label{maincostum2}	
Let $V^n$ and $U^p$ denote real vector 
spaces where $V^n$ carries $J\in End(V)$ satisfying $J^2=-I$ 
and $U^p$ is endowed with a positive definite inner product. 
Given a symmetric bilinear form \mbox{$\a\colon V^n\times V^n\to U^p$} 
let $\gamma\colon V^n\times V^n\to U^p$ be given by 
$$
\gamma(X,Y)=\alpha(X,Y)+\alpha(JX,JY).
$$
Assume that $\beta\colon V^n\times V^n\to W^{p,p}=U^p\oplus U^p$ 
defined by 
\be\label{betagamma}
\beta(X,Y)=(\gamma(X,Y),\gamma(X,JY))
\ee  
be flat with respect to the inner product on $W^{p,p}$ given by
$$
\lp(\xi_1,\xi_2),(\eta_1,\eta_2)\rp
=\<\xi_1,\eta_1\>_{U^p}-\<\xi_2,\eta_2\>_{U^p}.
$$
Then $\dim\mathcal{N}(\beta)=n-r$ where $r$ was defined by \eqref{r}.
\end{proposition}

\proof First observe that the symmetric bilinear form 
$\gamma$ satisfies
$$
\gamma(X,JY)+\gamma(JX,Y)=0
$$
for any $X,Y\in V^n$. Then 
$$
\lp\beta(X,Y),\beta(Y,X)\rp=\|\gamma(X,Y)\|^2+\|\gamma(X,JY)\|^2.
$$
Thus
\be\label{e:positivedefi}
\beta(X,Y)=0\;\;\mbox{and}\;\;\text{if and only if}
\;\;\lp\beta(X,Y),\beta(Y,X)\rp=0
\ee
for $X,Y\in V^n$.

Fix $X\in RE(\beta)$ and set $N=\ker B_X$. Since
$\mathcal{S}(\beta|_{N\times N})\subset\U(X)$ by Lemma \ref{bilfla}, 
then 
$$
\lp\beta(\eta_1,\eta_2),\beta(\eta_2,\eta_1)\rp=0
$$
for any $\eta_1,\eta_2\in N$. It follows from \eqref{e:positivedefi} 
that $\beta|_{N\times N}=0$. Now flatness of $\beta$ yields
$$
\lp\beta(Y,\eta),\beta(\eta,Y)\rp
=\lp\beta(Y,Y),\beta(\eta,\eta)\rp=0
$$
for any $\eta\in N$ and $Y\in V$. 
Then \eqref{e:positivedefi} gives $\beta(Y,\eta)=0$, and  
hence $N=\mathcal{N}(\beta)$. We have  
$$
\dim\mathcal{N}(\beta)=\dim N=n-\dim B_X(V)=n-r,
$$
as we wished.\vspace{2ex}\qed

\noindent {\em Proof of Theorem \ref{minimal}}: Let  
$\beta\colon T_xM\times T_xM\to W^{p,p}=N_fM(x)\oplus N_fM(x)$ 
be given by 
$$
\beta(X,Y)=(\a(X,Y),\a(X,JY)),
$$
where $\a\colon T_xM\times T_xM\to N_fM(x)$ denotes the second 
fundamental form of $f$ at $x\in M^{2n}$ and $W^{p,p}$ is endowed 
with the inner product
$$
\lp(\xi_1,\xi_2),(\eta_1,\eta_2)\rp
=\<\xi_1,\eta_1\>_{N_fM(x)}-\<\xi_2,\eta_2\>_{N_fM(x)}.
$$
Let $\gamma,\tilde{\gamma}\colon T_xM\times T_xM\to N_fM(x)$
be the symmetric bilinear forms defined by
$$
\gamma(X,Y)=\alpha(X,Y)+\alpha(JX,JY)\;\;\mbox{and}
\;\;\tilde{\gamma}(X,Y)=\alpha(X,Y)-\alpha(JX,JY) 
$$
and let $\beta^\gamma,\beta^{\tilde{\gamma}}
\colon T_xM\times T_xM\to W^{p,p}$ be the bilinear forms 
associated by \eqref{betagamma} to $\gamma,\tilde{\gamma}$,
respectively, that is,
$$
\beta^\gamma(X,Y)=(\gamma(X,Y),\gamma(X,JY))
\;\;\mbox{and}\;\;\beta^{\tilde\gamma}(X,Y)
=(\tilde\gamma(X,Y),\tilde\gamma(X,JY)).
$$
Then
\be\label{e:betasum}
2\beta=\beta^\gamma+\beta^{\tilde{\gamma}}.
\ee
To conclude the proof, it suffices to show that $\beta^\gamma=0$
since this is equivalent to
\be\label{pluri}
\alpha(X,JY)=\alpha(Y,JX)\;\;\mbox{for all}\;\;X,Y\in T_xM
\ee 
and, in particular, to $f$ being minimal.

Using the Gauss equation for $f$ and that the curvature tensor 
of $M^{2n}$ satisfies
$$
J\circ R(X,Y)=R(X,Y)\circ J\;\;\text{and}\;\;R(JX,JY)=R(X,Y), 
$$ 
it is easy to verify that $\beta^\gamma$ is a flat bilinear form 
and that
\be\label{comb}
\lp\beta^\gamma(X,Y),\beta^{\tilde{\gamma}}(Z,T)\rp
=\lp\beta^\gamma(X,T),\beta^{\tilde{\gamma}}(Z,Y)\rp
\ee
for all $X,Y,Z,T\in T_xM$.  

If $(\xi,\bar{\xi})\in\mathcal{S}(\beta^\gamma)$, 
then also $(\xi,-\bar{\xi}),(\bar{\xi},\xi),
(\bar{\xi},-\xi)\in\mathcal{S}(\beta^\gamma)$. In fact, let
$$
(\xi,\bar{\xi})=\sum_k\beta^\gamma(X_k,Y_k)
=\sum_k(\alpha(X_k,Y_k)+\alpha(JX_k,JY_k),
\alpha(X_k,JY_k)-\alpha(JX_k,Y_k)).
$$
Then
$$
\sum_k\beta^\gamma(Y_k,X_k)=(\xi,-\bar{\xi}),\;\;
\sum_k\beta^\gamma(JY_k,X_k)=(\bar{\xi},\xi)\;\text{and}\;
\sum_k\beta^\gamma(X_k,JY_k)=(\bar{\xi},-\xi).
$$
We have just shown that
$$
\pi_1(\mathcal{S}(\beta^\gamma))
=\pi_2(\mathcal{S}(\beta^\gamma)),
$$
where $\pi_j\colon W^{p,p}\to N_fM(x)$, $j=1,2$, denote 
the projections onto its components.

Setting $U^\gamma=\pi_1(\mathcal{S}(\beta^\gamma))$, we 
claim that
\be\label{e:imagebetagamma}
\mathcal{S}(\beta^\gamma)=U^\gamma\oplus U^\gamma.
\ee
In fact, we know that 
$\mathcal{S}(\beta^\gamma)\subset U^\gamma\oplus U^\gamma$. 
On the other hand, let $(\xi,\bar{\eta})\in U^\gamma\oplus U^\gamma$. 
By the definition of $U^\gamma$, there exist $\bar{\xi},\eta\in U^\gamma$ 
such that $(\xi,\bar{\xi}),(\eta,\bar{\eta})\in\mathcal{S}(\beta^\gamma)$.  
But then  $(\xi,-\bar{\xi})$ and $(-\eta,\bar{\eta})$ also
belong to $\mathcal{S}(\beta^\gamma)$, and thus 
$(\xi,\bar{\eta})\in\mathcal{S}(\beta^\gamma)$.  This proves the claim.

Assume that $\dim U^\gamma=s>0$. It follows 
from \eqref{e:betasum} that
\be\label{betas}
2\beta|_{V\times\mathcal{N}(\beta^\gamma)}
=\beta^{\tilde{\gamma}}|_{V\times\mathcal{N}(\beta^\gamma)}
\ee
and from \eqref{comb} that
\be\label{imagebe}
\spa\{\beta^{\tilde{\gamma}}(X,Z): X\in V^n\;
\text{and}\;Z\in\mathcal{N}(\beta^\gamma)\}
\subset\mathcal{S}(\beta^\gamma)^\perp.
\ee
Let $\{\xi_1,\dots,\xi_s\}$ be an orthonormal basis of $U^\gamma$.
We obtain from \eqref{e:imagebetagamma}, \eqref{betas} and 
\eqref{imagebe} that
$$
\<\alpha(X,Y),\xi_j\>=\lp\beta(X,Y),(\xi_j,0)\rp=0
$$
and 
$$
\<\alpha(X,JY),\xi_j\>=-\lp\beta(X,Y),(0,\xi_j)\rp=0
$$
for all $X\in T_xM$ and $Y\in\mathcal{N}(\beta^\gamma)$.  
Therefore 
$$
\alpha_{U^\gamma}(X,Y)=0=\alpha_{U^\gamma}(X,JY)
$$
for all $X\in T_xM$ and $Y\in\mathcal{N}(\beta^\gamma)$. 
Hence $\nu_s^c\geq\dim\mathcal{N}(\beta^\gamma)$.  Since
Proposition \ref{maincostum2} gives
$$
\dim\mathcal{N}(\beta^\gamma)\geq 2n-\dim\mathcal{S}(\beta^\gamma),
$$ 
thus $\nu_s^c\geq2(n-s)$. This contradicts the assumptions on the 
complex $s$-nullities and proves that $U^\gamma=0$, and hence
$\beta^\gamma=0$.\qed

\begin{remarks}\po {\em  
$(1)$ In order to conclude minimality it is not 
unexpected to prove \eqref{pluri} which seems a quite stronger 
condition. In fact, it was shown in \cite{DR1} that \eqref{pluri} 
and minimality are equivalent conditions.  \vspace{1ex}

\noindent $(2)$ It was shown in \cite{DF} that any substantial
real Kaehler submanifold $f\colon M^{2n}\to\R^{2n+p}$, $p\geq 2$, 
without local Euclidean factor, free of flat points and such
that $\nu^c_p=2n-2$ has to be minimal.
}\end{remarks}

\begin{lemma}\po\label{l:almostcomplex}
Let $V^n$ be a real vector space that carry $T\in End(V^n)$ satisfying 
$T^2=-I$. If $X_1,TX_1,\ldots,X_{k-1},TX_{k-1},X_k$ are linearly 
independent then $X_1,TX_1,\ldots,X_k,TX_k$ are also linearly independent.  
In particular $n$ is even.	
\end{lemma}

\proof If otherwise, there is 
$0\neq(a_1,\ldots,a_k,b_1,\ldots,b_{k-1})\in\R^{2k-1}$ such that
$$
TX_k=\sum_{i=1}^ka_iX_i+\sum_{j=1}^{k-1}b_jTX_j.
$$
It follows easily that
$$
(1+a_k^2)X_k + \sum_{j=1}^{k-1}((a_ka_j-b_j)X_j+(a_kb_j+a_j)TX_j)=0,
$$
and this is a contradiction.\qed 

\begin{proposition}\po\label{maincostum}
Let $V^n$ and $U^p$ be real vector 
spaces where $V^n$ carries \mbox{$J\in End(V^n)$} satisfying $J^2=-I$ 
and $U^p$ is endowed with a positive definite inner product.
Then let $\a\colon V^n\times V^n\to U^p$, $n\geq 2p$, be a symmetric 
bilinear form satisfying 
$$
\a(JX,Y)=\a(X,JY)
$$
for $X,Y\in V^n$. Assume that the symmetric bilinear form
$\beta\colon V^n\times V^n\to W^{p,p}$ defined~by 
\be\label{beta}
\beta(X,Y)=(\a(X,Y),\a(X,JY))
\ee  
is flat with respect to the inner product on $W^{p,p}=U^p\oplus U^p$ 
given by
$$
\lp(\xi_1,\xi_2),(\eta_1,\eta_2)\rp
=\<\xi_1,\eta_1\>_{U^p}-\<\xi_2,\eta_2\>_{U^p}.
$$
If $p\leq 11$ and the subspace $\mathcal{S}(\beta)$ is nondegenerate 
then $\dim\mathcal{N}(\beta)\geq n-2p$.
\end{proposition}

\proof 
Fix $X\in RE^o(\beta)$ and denote $N(X)=\ker B_X$.  We first argue 
that $\tau\leq p-1$. If otherwise, the density of $RE^o(\beta)$
gives
$$
\lp\beta(Z,Y),\beta(Z,T)\rp=0
$$
for any $Z,Y,T\in V^n$. Then flatness yields
$$
0=\lp\beta(Z+R,Y),\beta(Z+R,T)\rp= 2\lp\beta(Z,Y),\beta(R,T)\rp
$$
for any $Z,Y,T,R\in V^n$, and that contradicts that the subspace
$\mathcal{S}(\beta)$ is nondegenerate.  

If $\tau=0$ then $\dim\mathcal{N}(\beta)\geq n-2p$ holds. In fact, 
now Lemma \ref{bilfla} yields $N(X)=\mathcal{N}(\beta)$ and
therefore $\dim\mathcal{N}(\beta)=\dim N(X)\geq n-2p$.  
Consequently, in the sequel we work with $1\leq\tau\leq p-1$. 

We have that  
\be\label{iff}
(\eta,\bar{\eta})\in\U(X)\;\;\text{if and only if}\;\; 
(\bar{\eta},-\eta)\in\U(X).
\ee
In fact, if $(\eta,\bar{\eta})=B_XZ$ then 
$(\bar{\eta},-\eta)=B_XJZ$. Moreover, since 
$$
0=\lp B_XY,(\eta,\bar{\eta})\rp
=\<\a(X,Y),\eta\>-\<\a(X,JY),\bar{\eta}\>
$$ 
for any $Y\in V^n$, then also
$$
\lp B_XY,(\bar{\eta},-\eta)\rp
=\<\a(X,Y),\bar{\eta}\>+\<\a(X,JY),\eta\>=0
$$
for any $Y\in V^n$. 

We have from \eqref{iff} that $T\in End(\U(X))$ given by 
$T(\eta,\bar{\eta})=(\bar{\eta},-\eta)$ is well defined. 
Since $T^2=-I$, it follows from Lemma \ref{l:almostcomplex} that 
$\tau=2s$ and that we may write
\be\label{base}
\U(X)=\spa\{(\eta_1,\bar{\eta}_1),(\bar{\eta}_1,-\eta_1), 
\dots,(\eta_s,\bar{\eta}_s),(\bar{\eta}_s,-\eta_s)\}.
\ee   

Decomposing $W^{p,p}$ as in \eqref{e:decomposition}, we have 
$$
W^{p,p}=\U(X)\oplus\hat{\U}(X)\oplus\mathcal{V}^{p-\tau,p-\tau}
$$
and $B_X(V^n)\subset\U(X)\oplus\mathcal{V}$. Let 
$\hat{\beta}$ denote the $\hat{\U}(X)$-component of $\beta$. 
That the subspace $\mathcal{S}(\beta)$ is nondegenerate gives 
$\mathcal{S}(\hat{\beta})=\hat{\U}(X)$. 

Setting $\hat{B}_YZ=\hat{\beta}(Y,Z)$, it is easy to see that
\be\label{maisys}
\begin{cases}
\lp\hat{B}_{Y}Z,(\eta_j,\bar{\eta}_j)\rp 
=-\lp\hat{B}_{Y}JZ,(\bar{\eta}_j,-\eta_j)\rp,\;1\leq j\leq s,
\vspace{1ex}\\
\lp\hat{B}_{Y}Z,(\bar{\eta}_j,-\eta_j)\rp 
=\lp\hat{B}_{Y}JZ,(\eta_j,\bar{\eta}_j)\rp,\;1\leq j\leq s,
\end{cases}
\ee 
for any $Y,Z\in V^n$. It follows from \eqref{base} and \eqref{maisys} 
that $S\in End(\hat{\U}(X))$ given by $S\hat{B}_Y(Z)=\hat{B}_Y(JZ)$ 
is well defined. Since $S^2=-I$, then Lemma \ref{l:almostcomplex} 
gives the following:

\begin{fact}\po\label{fact} If the vectors 
$\hat{B}_{Y_1}Z_1,\hat{B}_{Y_1}JZ_1,\ldots,\hat{B}_{Y_r}Z_r$ 
are linearly independent, then the same is true for 
$\hat{B}_{Y_1}Z_1,\hat{B}_{Y_1}JZ_1,\ldots,
\hat{B}_{Y_r}Z_r,\hat{B}_{Y_r}JZ_r$. 
\end{fact}

Take $Y\in RE(\hat{\beta})$. Then Fact \ref{fact} gives that 
$\kappa=\dim\hat{B}_Y(V^n)$ is even. We show next that if 
$\kappa=\tau$, then $\dim\mathcal{N}(\beta)\geq n-2p$ 
holds regardless the value of $p$.  We have 
$\hat{B}_Y(V^n)=\hat{\U}(X)$.  
Set $B_1=B_Y|_{N(X)}\colon N(X)\to\U(X)$ and $N_1=\ker B_1$.   
Then $\dim N(X)\leq\dim N_1+\tau$. By Lemma \ref{bilfla}, 
$$
B_Z\eta=0\;\;\text{if and only if}\;\;\lp B_Z\eta,\hat{B}_Y(V^n)\rp=0
$$
for $\eta\in N_1$ and $Z\in V^n$. Flatness of $\beta$ 
gives
$$
\lp B_Z\eta,\hat{B}_Y(V^n)\rp=\lp B_Z\eta,B_Y(V^n)\rp
=\lp B_Z(V^n),B_Y\eta\rp=0
$$
for $\eta\in N_1$ and any $Z\in V^n$.
Thus $N_1\subset\mathcal{N}(\beta)$, and hence
$$
\dim\mathcal{N}(\beta)\geq\dim N_1\geq\dim N(X)
-\tau\geq n-(2p-\tau)-\tau=n-2p.
$$

On one hand, as seen above we only have to consider the cases 
$2\leq\kappa<\tau\leq p-1$ where $\kappa$ and $\tau$ are even. On the 
other hand, being $\beta$ symmetric then Lemma \ref{l:symmetricspan} 
applied to $\hat{\beta}$ gives $\kappa(\kappa+1)\geq 2\tau$. In particular, 
it suffices to argue for $\tau=6,8,10$.  

\begin{fact}\po\label{fevendim}
If  $L\subset V^n$ is a $J$-invariant subspace then  
the subspace $B_Y(L)$ has even dimension for any $Y\in V^n$.	
\end{fact}

In fact, we have that $T_0\in End(B_Y(L))$ 
given by $T_0B_YZ=B_YJZ$ for $Z\in L$ is well defined and satisfies 
$T_0^2=-I$, and the claim follows from Lemma \ref{l:almostcomplex}.
\vspace{1ex}

\noindent \emph{Case $\tau=6$ and $\kappa=4$}. By Fact \ref{fact} 
there are vectors $Y_1,Y_2$ in the open dense subset 
$RE^o(\beta)\cap RE(\hat{\beta})$ of $V^n$ such that 
\be\label{e:spanuchapel}
\hat{\U}(X)=\hat{B}_{Y_1}(V^n)\oplus\hat{B}_{Y_2}(V^n).
\ee
First suppose that $Y_1$ satisfies $\dim B_{Y_1}(N(X))\leq 4$.  
Set $B_1=B_{Y_1}|_{N(X)} \colon N(X)\to\U(X)$ and $N_1=\ker B_1$. 
Then $\dim N_1\geq\dim N(X)-4$. Flatness and \eqref{bili} easily 
yield 
$$
\rank B_{Y_2}(N_1)\oplus\hat{B}_{Y_1}(V^n) = 0.
$$
Thus $\dim B_{Y_2}(N_1)\leq 2$.  
Set $B_2=B_{Y_2}|_{N_1}\colon N_1\to\U(X)$ and $N_2=\ker B_2$.  
Then $\dim N_2\geq\dim N_1-2$.  
We claim that $N_2\subset\mathcal{N}(\beta)$. 
Now \eqref{bili} and \eqref{e:spanuchapel} give
$$
B_Z\eta=0\;\;\text{if and only if}\;\;
\lp B_Z\eta,\hat{B}_{Y_j}(V^n)\rp=0,\;j=1,2,
$$
for any $\eta\in N_2$. 
Flatness of $\beta$ and \eqref{bili} yield
$$
\lp B_Z\eta,\hat{B}_{Y_j}(V^n)\rp=\lp\beta(Z,V^n),\beta(Y_j,\eta)\rp=0
$$
for any $Z\in V^n$, $j=1,2$, and the claim has been proved. Then
$$
\dim\mathcal{N}(\beta)\geq\dim N_2\geq\dim N_1-2\geq\dim N(X)-6\geq n-2p.
$$

To conclude the proof of this case, in view of Fact \ref{fevendim}
it remains to show that we cannot have $B_Y(N(X))=\U(X)$ for any 
$X\in RE^o(\beta)$ and $Y\in RE^o(\beta)\cap RE(\hat{\beta})$.  
Assume otherwise, and consider the decomposition 
$$
W^{p,p}=\U(Y_1)\oplus\hat{\U}(Y_1)\oplus\mathcal{V}(Y_1)
$$
given by \eqref{e:decomposition}. Let $\tilde{\beta}$ be the 
$\hat{\U}(Y_1)$ component of $\beta$ and 
$Y_2\in RE^o(\beta)\cap RE(\hat{\beta})\cap RE(\tilde{\beta})$ be
such that \eqref{e:spanuchapel} holds. Then $B_{Y_2}(N(Y_1))=\U(Y_1)$ 
and hence
$$
B_{Y_1}(V^n)\cap\U(X)=\U(X),\;\;
B_{Y_2}(V^n)\cap\U(X)=\U(X)\;\;\mbox{and}\;\;
B_{Y_2}(V^n)\cap\U(Y_1)=\U(Y_1).
$$
It follows that $\U(X)+\U(Y_1)+\U(Y_2)$ is an isotropic subspace.  
Moreover, since we have that $\dim\U(X)\cap\U(Y_j)=\dim\U(Y_1)\cap\U(Y_2)=2$,
then its dimension would be $12$, and this is  a contradiction because 
$p\leq 11$.  
\vspace{1ex}

\noindent \emph{Case $\tau=8$ and $\kappa=4$}. 
Let $Y_1\in RE^o(\beta)\cap RE(\hat{\beta})$ and decompose
$\hat{\U}=\hat{B}_{Y_1}(V^n)\oplus\hat{\U}_1$.
Let $\beta^1$ be the $\hat{\U}_1$-component of $\beta$ and set 
$B^1_Y(Z)=\beta^1(Y,Z)$. Since $\mathcal{S}(\beta^1)=\hat{\U}_1$,  
we have from Lemma \ref{l:symmetricspan} that $3\leq\dim B_{Y}^1(V^n)\leq 4$ 
for $Y\in RE(\beta^1)$.  We argue that $\dim B_{Y}^1(V^n)=4$. 
If otherwise, then 
$$
\dim\hat{B}_{Y_1}(V^n)+\hat{B}_{Y_2}(V^n)=7
$$
for $Y_2\in RE^o(\beta)\cap RE(\hat{\beta})\cap RE(\beta^1)$.  
Thus $\dim\hat{B}_{Y_1}(V^n)\cap\hat{B}_{Y_2}(V^n)=1$ since $\kappa=4$.  
To see that this cannot happen, suppose $\hat{B}_{Y_1}(X_1)=\hat{B}_{Y_2}(X_2)$. 
We have from \eqref{maisys} that $\hat{B}_{Y_1}(JX_1)=\hat{B}_{Y_2}(JX_2)$.  
Therefore $\hat{B}_{Y_1}(X_1)=a\hat{B}_{Y_1}(JX_1)$ and \eqref{maisys} 
gives $\hat{B}_{Y_1}(X_1)=0$, a contradiction.  We conclude that there exist 
$Y_1,Y_2\in RE^o(\beta)\cap RE(\hat{\beta})$ such that 
\eqref{e:spanuchapel} holds.

First suppose that $\dim B_{Y_1}(N(X))\leq 4$. 
Set $B_1=B_{Y_1}|_{N(X)}\colon N(X)\to\U(X)$ and 
$N_1=\ker B_1$.  Thus $\dim N(X)\leq\dim N_1+4$.  
Set $B_2=B_{Y_2}|_{N_1}\colon N_1\to\U(X)$ and $N_2=\ker B_2$.  
 From flatness and \eqref{bili} we have
$$
\lp B_{Y_2}\eta,\hat{B}_{Y_1}(V^n)\rp=\lp B_{Y_2}(V^n),B_{Y_1}\eta\rp=0
$$
for any $\eta\in N_1$, and thus $\dim B_{Y_2}(N_1)\leq 4$.  
Hence, $\dim N_1\leq\dim N_2+4$.  We claim that 
$N_2\subset\mathcal{N}(\beta)$. 
Now \eqref{bili} and \eqref{e:spanuchapel} give
$$
B_Z\eta=0\;\;\text{if and only if}\;\;
\lp B_Z\eta,\hat{B}_{Y_j}(V^n)\rp=0,\;j=1,2,
$$
for any $\eta\in N_2$. Flatness of $\beta$ and \eqref{bili} give
$$
\lp B_Z\eta,\hat{B}_{Y_j}(V^n)\rp=\lp\beta(Z,V^n),\beta(Y_j,\eta)\rp=0
$$
for any $Z\in V^n$, $j=1,2$, and the claim has been proved.  Then
$$
\dim\mathcal{N}(\beta)\geq\dim N_2\geq\dim N_1-4\geq\dim N(X)-8\geq n-2p.
$$

Observe that we cannot have $B_{Y_1}(N(X))=\U(X)$ since being
$\kappa=4$ and $p\leq 11$, we would have that $\dim\U(Y_1)\leq 7$.
Hence, by Fact \ref{fevendim} the remaining case to be examined 
is when $\dim B_{Y_1}(N(X))=6$. Set 
$B_1=B_{Y_1}|_{N(X)}\colon N(X)\to\U(X)$ and $N_1=\ker B_1$. 
Flatness of $\beta$ and \eqref{bili} give
$$
\lp B_{Y_2}(N_1),\hat{B}_{Y_1}(V^n)\rp=0
$$
and thus $\dim B_{Y_2}(N_1)\leq 4$. But in fact, 
$\dim B_{Y_2}(N_1)\leq 2$ since otherwise $\dim\U(Y_2)\leq 7$.   
Set $B_2=B_{Y_2}|_{N_1}\colon N_1\to\U(X)$ 
and $N_2=\ker B_2$.  As above, 
$N_2\subset\mathcal{N}(\beta)$ and  then 
$$
\dim\mathcal{N}(\beta)\geq\dim N_2\geq\dim N_1-2\geq\dim N(X)-8\geq n-2p.
$$

\noindent \emph{Case $\tau=8$ and $\kappa=6$}. 
Let $Y_1,Y_2\in RE^o(\beta)\cap RE(\hat{\beta})$ be such that 
\eqref{e:spanuchapel} is satisfied. Then $\dim B_{Y_j}(N(X))\leq 4$ 
for $j=1,2$.  Otherwise, if $\dim B_{Y_1}(N(X))\geq 6$ then 
$\dim\U(Y_1)\leq 7$, a contradiction since $p\leq 11$. Again let 
$B_1=B_{Y_1}|_{N(X)}\colon N(X)\to\U(X)$, 
$N_1=\ker B_1$, $B_2=B_{Y_2}|_{N_1}\colon N_1\to\U(X)$ and 
$N_2=\ker B_2$. As above, we have $N_2\subset\mathcal{N}(\beta)$. 
Hence 
$$
\dim\mathcal{N}(\beta)\geq\dim N_2\geq\dim N_1-4
\geq\dim N(X)-8\geq n-2p.
$$

\noindent \emph{Case $\tau=10$}. 
If $B_Y(Z)\in\U(X)$ then \eqref{iff} gives that also $B_YJZ\in\U(X)$. Moreover, 
computing the inner products shows that
$$
\lp B_YZ,\hat{B}_RT\rp=-\lp B_YJZ,\hat{B}_RJT\rp
$$
for any $R,T\in V^n$.  Then Lemma \ref{l:almostcomplex} and Fact 
\ref{fevendim} give
$$
\rank B_Y(V^n)\cap\U(X)\oplus\hat{B}_Y(V^n)=4\rho,
$$
where $\rho\in\mathbb{N}$. But since $\tau=10$ and 
hence $p=11$, then $\rho=0$, that is,
\be\label{rankestimate10}
\rank B_Y(V^n)\cap\U(X)\oplus\hat{B}_Y(V^n)=0
\ee
for any $Y\in RE^o(\beta)\cap RE(\hat{\beta})$.  We argue for $\kappa=4$ 
being the other cases similar. 
An analogous argument as the one given for the case $\tau=8$ and $\kappa=4$ 
shows that there is $Y_1,Y_2,Y_3\in RE^o(\beta)\cap RE(\hat{\beta})$ 
such that
\be\label{spanuchapel3}
\hat{\U}(X)=\hat{B}_{Y_1}(V^n)+\hat{B}_{Y_2}(V^n)+\hat{B}_{Y_3}(V^n)
\ee
where 
$$
\dim\hat{B}_{Y_1}(V^n)+\hat{B}_{Y_2}(V^n)=8.
$$

Set $B_1=B_{Y_1}|_{N(X)}\colon N(X)\to\U(X)$, 
$N_1=\ker B_1$, \mbox{$B_2=B_{Y_2}|_{N_1}\colon N_1\to\U(X)$}, 
$N_2=\ker B_2$ and $B_3=B_{Y_3}|_{N_2}\colon N_2\to\U(X)$. 
From \eqref{rankestimate10} we have $\dim B_{Y_1}(N(X))\leq 6$.  
Flatness of $\beta$ and \eqref{bili} give
$$
\lp B_{Y_2}\eta,\hat{B}_{Y_1}(V^n)\rp=\lp B_{Y_2}(V^n),B_{Y_1}\eta\rp=0
$$ 
for any $\eta\in N_1$.  Thus from \eqref{rankestimate10} we 
have $\dim B_{Y_2}(N_1)\leq2$.  Using again flatness of $\beta$ and 
\eqref{bili} we obtain
$$
\lp B_{Y_3}(N_2),\hat{B}_{Y_j}(V^n)\rp=0,\;\;\mbox{j=1,2}.
$$
Hence \eqref{rankestimate10} gives $B_{Y_3}(N_2)=0$.  
As above, we have $N_2\subset\mathcal{N}(\beta)$.  Then
$$
\dim\mathcal{N}(\beta)\geq\dim N_2\geq\dim N_1-2
\geq \dim N(X)-8\geq n-2p, 
$$
and this concludes the proof.\qed

\begin{corollary}\po\label{mainlemma} Assume that the bilinear form 
$\beta\colon V^n\times V^n\to W^{p,p}$ given by \eqref{beta} 
is flat. If $p\leq 11$ and $\dim\mathcal{N}(\beta)\leq n-2p-1$, then 
$\U=\mathcal{S}(\beta)\cap\mathcal{S}(\beta)^\perp$ satisfies
$\dim\U=s>0$ is even.  Moreover, the projection $U_1$ of $\U$ 
on the first component of $W^{p,p}$ satisfies $\dim U_1=s$ and the 
$U_2$ component $\a_2$ of $\alpha$ with respect to the orthogonal 
decomposition $U^p=U_1^s\oplus U_2^{p-s}$ satisfies
$$
\dim\mathcal{N}(\a_2)\cap J\mathcal{N}(\a_2)\geq n-2(p-s).
$$
\end{corollary}

\proof Proposition \ref{maincostum} gives $s>0$. 
If $(\xi,\bar{\xi})\in\U$ then
$$
(\xi,\bar{\xi})=\sum_i\beta(X_i,Y_i)
=\sum_i(\a(X_i,Y_i),\a(X_i,JY_i))
$$
and
$$
\lp\beta(X,Y),(\xi,\bar{\xi})\rp=\<\a(X,Y),\xi\>
-\<\a(X,JY),\bar{\xi}\>=0
$$
for all $X,Y\in V^n$.  Hence 
$(\bar{\xi},-\xi)=\sum_i\beta(X_i,JY_i)\in\mathcal{S}(\beta)$
and 
$$
\lp(\beta(X,Y),(\bar{\xi},-\xi)\rp=\<\a(X,Y),\bar{\xi}\>
+\<\a(X,JY),\xi\>=0
$$
for all $X,Y\in V^n$. Thus $(\bar{\xi},-\xi)\in\U$, and hence 
it follows easily from of Lemma \ref{l:almostcomplex} that 
$s$ is even. Moreover, we have $\pi_1(\U)=U_1^s=\pi_2(\U)$
where $\pi_j\colon W^{p,p}\to U^p$, $j=1,2$, denote the 
projections.

Set $\beta=\beta_1+\beta_2$ where 
$\beta_j\colon V^n\times V^n\to U_j\oplus U_j$ is given by 
$$
\beta_j(X,Y)=(\a_j(X,Y),\a_j(X,JY))
$$
and $\a_1=\pi_{U_1}\circ\a$. Then $\mathcal{S}(\beta_1)=\U$. Since 
$\beta$ is flat and $\beta_1$ is null then $\beta_2$ is flat. 
Moreover, since the subspace $\mathcal{S}(\beta_2)$ is
nondegenerate,  from Proposition \ref{maincostum} we obtain
$$
\dim\mathcal{N}(\beta_2)\geq n-2(p-s).
$$
To conclude the proof observe that
$\mathcal{N}(\beta_2)=\mathcal{N}(\a_2)\cap J\mathcal{N}(\a_2)$.
\vspace{2ex}\qed

\noindent {\em Proof of Theorem \ref{main}:}
We make use of some of the contents of the proof of Theorem~\ref{minimal}. 
We have seen that $\beta$ is symmetric. Moreover, a similar argument 
as the one used in that proof gives that $\beta$ is flat.
We have from \eqref{pluri} that the complex $s$-nullities 
coincide with the standard $s$-nullities defined in \cite{CD}, that is, 
$\nu^c_s(x)=\nu_s(x)$ where
$$
\nu_s(x)=\max_{U^s\subset N_fM(x)}\{\dim\mathcal{N}(\a_{U^s})\}.
$$
In particular $\dim\mathcal{N}(\beta)\leq 2(n-p)-1$, and hence
Corollary \ref{mainlemma} applies. Therefore, 
$\dim\U=s>0$, and there is an orthogonal decomposition
$N_fM(x)=U_1^s\oplus U_2^{p-s}$ such that
$$
\dim\mathcal{N}(\alpha_2)\cap J\mathcal{N}(\alpha_2)\geq 2(n-p+s).
$$
Thus $\nu_{p-s}(x)\geq 2(n-p+s)$, and this is a contradiction 
unless $p=s$ and $U_2=0$, that is, the bilinear form 
$\beta$ is null.

We know from \cite{DG1} that any simply-connected minimal 
real Kaehler submanifold has a  one-parameter associated family 
of minimal isometric immersions 
$f_\theta\colon M^{2n}\to\R^{2n+p}$, $\theta\in [0,\pi)$, with second 
fundamental form
$$
\alpha_\theta(X,Y)=\cos\theta\alpha(X,Y)+\sin\theta\alpha(JX,Y).
$$
Moreover, that family is trivial, that is, any pair of elements in the
family are congruent in $\R^{2n+p}$, if and only if $f$ is holomorphic. 
On the other hand, since $\nu_1(x)\leq 2n-3$ then
$\mathcal{S}(\alpha)=N_fM(x)$. Hence, being
$\beta$ null there is a smooth vector bundle isometry 
$T\colon N_fM\to N_{f_{\pi/2}}M$ such that
$T\circ\alpha=\alpha_{\pi/2}$. Using Theorem $3$ in \cite{No} 
or Lemma $4.16$ in \cite{DT} we have that $T$ is parallel in the
normal connection, and hence $f$ and $f_{\pi/2}$ are 
congruent. Thus, the one parameter associated family of minimal 
isometric immersions is trivial; see Proposition $1.6$ in \cite{DG1} 
or Theorem $15.10$ in \cite{DT}.\qed

\begin{remarks}\label{remark}\po {\em 
$(1)$  For codimension $p\leq 5$, the proof of Theorem \ref{main}
follows immediately from Theorem \ref{minimal}  and the rigidity 
result given in \cite{CD} since these results assure that the 
associated family discussed above is trivial.\vspace{1ex}

\noindent $(2)$ We point out that Theorem \ref{main} holds for
codimension $p=4$ under the assumption $\nu^c_s(x)<2(n-s)$ just 
for $s=2,4$.  This can be proved using Proposition $6$ of
\cite{CCD} instead of Proposition \ref{maincostum}.
\vspace{1ex}

\noindent $(3)$ We observe that Lemma $3$ in \cite{DG4} may require 
a slightly stronger assumption than the one asked there in order to 
hold. On the other hand, to prove the theorem given there one
can replace the Lemma $3$ there by Corollary \ref{mainlemma} in this 
paper.
}\end{remarks}

\section{Counterexamples}

In this section, we show that Proposition \ref{maincostum} 
does not hold for $p\geq 12$.  This is achieved by constructing 
explicit counterexamples to that result. 
\vspace{2ex}

Let $\gamma\colon V_0^m\times V_0^m\to W_0^{q,q}$, 
$m>2q$, $q\geq 6$, be a flat symmetric bilinear form such that 
$\mathcal{S}(\gamma)=W_0^{q,q}$ and $\dim\mathcal{N}(\gamma)<m-2q$.   
We point out that explicit examples of such $\gamma$'s have been 
constructed in the proof of Theorem 2 in \cite{DF}.
Then set $V^n=V_0^m\oplus V_0^m$ and $W^{p,p}=W_0^{q,q}\oplus W_0^{q,q}$, 
where the latter is endowed with the inner product
$$
\lp(\zeta^1,\zeta^2),(\eta^1,\eta^2)\rp
=\<\zeta^1,\eta^1\>_{W_0^{q,q}}-\<\zeta^2,\eta^2\>_{W_0^{q,q}}.
$$

\begin{proposition}\label{p:beta}\po 
Let $\beta\colon V^n\times V^n\to W^{p,p}$ be the symmetric 
bilinear form given by
$$
\beta((X,Y),(Z,T))
=(\gamma(X,Z)-\gamma(Y,T),\gamma(X,T)+\gamma(Y,Z))
$$
for any $X,Y,Z,T\in V_0^m$.  
Then $\beta$ is flat with $\mathcal{S}(\beta)=W^{p,p}$ and  
$\dim\mathcal{N}(\beta)<n-2p$.
\end{proposition}

\proof The symmetry of $\beta$ follows from the symmetry of $\gamma$
and a straightforward computation gives that $\beta$ is flat.
To see that $\mathcal{S}(\beta)=W^{p,p}$, let 
$(\zeta^1,\zeta^2)\in W^{p,p}$ with $\zeta^1,\zeta^2\in W_0^{q,q}$.  
By assumption 
$$
\zeta^1=\sum_k\gamma(X_k,Y_k)\;\;\text{and}\;\;
\zeta^2=\sum_\ell\gamma(Z_\ell,T_\ell)
$$  
where $X_k,Y_k,Z_\ell,T_\ell\in V_0^m$. Then
$$
\sum_k\beta((X_k,0),(Y_k,0))+\sum_\ell\beta((Z_\ell,0),(0,T_\ell))
=(\zeta^1,\zeta^2).
$$
Finally, we have $(X,Y)\in\mathcal{N}(\beta)$ if and only if 
$$
\gamma(X,Z)-\gamma(Y,T)=0=\gamma(X,T)+\gamma(Y,Z)
$$
for any $Z,T\in V_0^m$. Clearly, the right-hand side is equivalent 
to $X,Y\in\mathcal{N}(\gamma)$.\vspace{2ex}\qed

To conclude the construction of a counterexample we show next
that $\beta$ given by Proposition \ref{p:beta} has the structure
required by Proposition \ref{maincostum}.
\vspace{1ex}

We fix an orthogonal decomposition $W^{q,q}_0=P_0\oplus N_0$
such that $\dim P_0=q$ and the induced inner product 
$\<\,,\,\>_{P_0}$  on $P_0$ is positive definite.  Hence 
the induced inner product $\<\,,\,\>_{N_0}$on $N_0$ is 
negative definite. Taking orthogonal components we 
have $\zeta=\zeta_{P_0}+\zeta_{N_0}$ if $\zeta\in W^{q,q}_0$.
Hence, we have an orthogonal decomposition 
$$
W^{p,p}=P\oplus N
$$
where $P=(P_0,0,0,N_0)$ and $N=(0,N_0,P_0,0)$. Hence, the induced 
inner product on $P$ is positive definite and negative definite 
on $N$. 

Let $\bar P=P_0\oplus N_0$ be endowed with the positive definite 
inner product
$$
\<(p_1,n_1),(p_2,n_2)\>_{\bar P}
=\<p_1,p_2\>_{P_0}-\<n_1,n_2\>_{N_0}
$$
where $p_1,p_2\in P_0$ and $n_1,n_2\in N_0$.  Then the map
$\phi\colon P\to\bar{P}$ defined by
$$
\phi(p_1,0,0,n_1)=(p_1,n_1)
$$
is an isometry.  

Endow the vector space $L^{p,p}=\bar P\oplus\bar P$ with the inner 
product $\<\,,\,\>_{L^{p,p}}$ of signature $(p,p)$ given by
$$
\<((p_1,n_1),(p_2,n_2)),((p_3,n_3),(p_4,n_4))\>_{L^{p,p}}
=\<(p_1,n_1),(p_3,n_3)\>_{\bar{P}}-\<(p_2,n_2),(p_4,n_4)\>_{\bar{P}}
$$
where $p_j\in P_0$ and $n_j\in N_0$, $j=1,2,3,4$.  
The map $\Tal\colon W^{p,p}\to L^{p,p}$ defined by
$$
\Tal(\zeta^1,\zeta^2)
=\Tal(\zeta^1_{P_0}+\zeta^1_{N_0},\zeta^2_{P_0}+\zeta^2_{N_0})
=((\zeta^1_{P_0},\zeta^2_{N_0}),(-\zeta^2_{P_0},\zeta^1_{N_0}))
$$
is an isometry. In fact, its inverse 
$\Tal^{-1}\colon L^{p,p}\to W^{p,p}$ is
$$
\Tal^{-1}((p_1,n_1),(p_2,n_2))=(p_1+n_2,-p_2+n_1)
$$
and
\begin{align*}
\<\Tal(\zeta^1,\zeta^2),\Tal(\eta^1,\eta^2)\>_{L^{p,p}}
&=\<((\zeta^1_{P_0},\zeta^2_{N_0}),(-\zeta^2_{P_0},\zeta^1_{N_0})),
((\eta^1_{P_0},\eta^2_{N_0}),(-\eta^2_{P_0},\eta^1_{N_0}))\>_{L^{p,p}}\\
&=\<\zeta^1_{P_0},\eta^1_{P_0}\>_{P_0}-\<\zeta^2_{N_0},\eta^2_{N_0}\>_{N_0}
-\<\zeta^2_{P_0},\eta^2_{P_0}\>_{P_0}+\<\zeta^1_{N_0},\eta^1_{N_0}\>_{N_0}\\
&=\lp(\zeta^1,\zeta^2),(\eta^1,\eta^2)\rp.
\end{align*}

Let $\a\colon V^n\times V^n\to\bar P$ be defined by 
$\a=\phi\circ\pi_P\circ\beta$, where $\pi_P\colon W^{p,p}\to P$ 
denotes the projection. Set $\gamma=\gamma_{P_0}+\gamma_{N_0}$ 
according to the  decomposition $W^{q,q}_0=P_0\oplus N_0$. Then
\begin{align*}
\alpha((X,Y),(Z,T))&=\phi\circ\pi_P\circ\beta((X,Y),(Z,T))\\
&=\phi\circ\pi_P(\gamma(X,Z)-\gamma(Y,T),\gamma(X,T)+\gamma(Y,Z))\\
&=(\gamma_{P_0}(X,Z)-\gamma_{P_0}(Y,T),\gamma_{N_0}(X,T)
+\gamma_{N_0}(Y,Z)).
\end{align*}

Let $J\in End(V)$ be given by $J(X,Y)=(-Y,X)$. Hence $J^2=-I$.  
Then
$$
\alpha((X,Y),J(Z,T))=\alpha(J(X,Y),(Z,T)),
$$
that is, $\alpha$ satisfies condition \eqref{pluri}.  
The bilinear form 
$\beta^\alpha\colon V^n\times V^n\to\bar P\oplus\bar P$ defined by 
\eqref{beta} is given by
\begin{align*}
\beta^\alpha((X,Y),(Z,T))
&=(\alpha((X,Y),(Z,T)),\alpha((X,Y),J(Z,T)))\\
&=((\gamma_{P_0}(X,Z)-\gamma_{P_0}(Y,T),
\gamma_{N_0}(X,T)+\gamma_{N_0}(Y,Z)),\\
&\quad\quad\quad\quad(-\gamma_{P_0}(X,T)-\gamma_{P_0}(Y,Z),
\gamma_{N_0}(X,Z)-\gamma_{N_0}(Y,T))).
\end{align*}
Then
$$
\Tal^{-1}\circ\beta^\alpha((X,Y),(Z,T))
=\beta((X,Y),(Z,T)).
$$
Thus $\beta^\alpha\colon V^n\times V^n\to L^{p,p}$ is a 
flat symmetric bilinear form such that 
$\mathcal{S}(\beta^\alpha)=L^{p,p}$ and satisfies
$\dim\mathcal{N}(\beta^\alpha)<n-2p$.

\noindent Alcides de Carvalho, Sergio Chion, Marcos Dajczer\\
IMPA -- Estrada Dona Castorina, 110\\
22460--320, Rio de Janeiro -- Brazil\\
e-mail: alcidesj@impa.br, sergio.chion@impa.br, marcos@impa.br
\end{document}